\documentclass[12pt]{article}
\usepackage{amsmath}
\usepackage{amssymb}
\usepackage{graphicx}

\def\be#1\ee{\begin{equation}#1\end{equation}}
\newcommand{\bea}{\begin{eqnarray}}
\newcommand{\eea}{\end{eqnarray}}
\newcommand{\beas}{\begin{eqnarray*}}
\newcommand{\eeas}{\end{eqnarray*}}

\newcommand{\hh}{\hat h}
\def\f{{\bf f}}
\def\n{{\bf n}}
\def\H{{\bf H}}
\def\J{{\bf J}}
\def\E{{\bf E}}
\def\r{{\bf r}}
\def\R{\mbox{\boldmath$\mathcal R$}}

\newtheorem{remark}{Remark}

\def\TA{\mbox{\boldmath$\mathcal A$}}

\newcommand{\hrarrow}  %horizontal right arrow 
{\linethickness{1mm}
\put(0,0){\line(1,0){3}}
\thicklines
\put(2,-1.5){\line(4,1){6}}
\put(2,1.5){\line(4,-1){6}}
\put(2,-1.5){\line(0,1){3}}
\put(2.4,-1.4){\line(0,1){2.8}}
\put(2.8,-1.3){\line(0,1){2.6}}
\put(3.2,-1.2){\line(0,1){2.4}}
\put(3.6,-1.1){\line(0,1){2.2}}
\put(4,-1.0){\line(0,1){2.0}}
\put(4.4,-0.9){\line(0,1){1.8}}
\put(4.8,-0.8){\line(0,1){1.6}}
\put(5.2,-0.7){\line(0,1){1.4}}
\put(5.6,-0.6){\line(0,1){1.2}}
\put(6,-0.5){\line(0,1){1.0}}
\put(6.4,-0.4){\line(0,1){0.8}}
\put(6.8,-0.3){\line(0,1){0.6}}
\thicklines}

\newcommand{\vdarrow} %vertical downward arrow
{\linethickness{1mm}
\put(0,0){\line(0,-1){3}}
\thicklines
\put(-1.5,-2){\line(1,-4){1.5}}
\put(1.5,-2){\line(-1,-4){1.5}}
\put(-1.5,-2){\line(1,0){3}}
\put(-1.4,-2.4){\line(1,0){2.8}}
\put(-1.3,-2.8){\line(1,0){2.6}}
\put(-1.2,-3.2){\line(1,0){2.4}}
\put(-1.1,-3.6){\line(1,0){2.2}}
\put(-1,-4){\line(1,0){2.0}}
\put(-0.9,-4.4){\line(1,0){1.8}}
\put(-0.8,-4.8){\line(1,0){1.6}}
\put(-0.7,-5.2){\line(1,0){1.4}}
\put(-0.6,-5.6){\line(1,0){1.2}}
\put(-0.5,-6){\line(1,0){1.0}}
\put(-0.4,-6.4){\line(1,0){0.8}}
\put(-0.3,-6.8){\line(1,0){0.6}}
\thicklines}

\begin{document}

\title{Fast finite-difference convolution for 3D problems in layered media}

\author{ Vladimir Druskin%
\thanks{Schlumberger Doll Research, 1 Hampshire St.,
Cambridge, MA 02139, USA (\texttt{druskin1@slb.com}).}
\and Mikhail Zaslavsky%
\thanks{Schlumberger Doll Research, 1 Hampshire St.,
Cambridge, MA 02139, USA (\texttt{mzaslavsky@slb.com}.)} }

\maketitle

\begin{abstract}
We developed fast direct solver for 3D Helmholtz and Maxwell equations in layered medium. The algorithm is based on the ideas of cyclic reduction for separable matrices. For the grids with major uniform part (within the survey domain in the problems of geophysical prospecting, for example) and small non-uniform part (PML and coarsening to approximate problems in infinite domain) the computational cost of our approach is $O(N_xN_ylog(N_xN_y)N_z)$. For general non-uniform grids the cost is $O(N^{3/2}_xN^{3/2}_yN_z)$. The first asymptotics coincide with the cost of FFT-based methods, which can be applied for uniform gridding (in $x$ and $y$) only. Our approach is significantly more efficient compared to the algorithms based on discrete Fourier transform which cost is $O(N^2_xN^2_yN_z)$. The algorithm can be easily extended for solving the elasticity problems as well. \end{abstract}
\section{Introduction}
With application to acoustic, elastic and electromagnetic scattering, the geological formations often can be approximated by  horizontally layered background media with several embedded inhomogeneities. For 3D case, such problems usually lead to the  large ill-conditioned system of linear equations when finite differences (or finite elements) are employed for discretization. In this work we assume that the discretization grid $S$ is a tensor product of one-dimensional grids $S_x\times S_y\times S_z$. Direct solvers can be efficiently applied for the solution of 2D (or 2.5D) problems as well as 3D problems of rather small size only. Therefore, the iterative solvers are the only option for 3D large-scale problems. The choice of preconditioner is crucial for robustness of iterative method. In \cite{zaslseg, zaslgeo} we show the efficiency of the preconditioning via discrete Green's function in layered medium. Similar ideas are used in the volume integral equation (IE) for the scattered field \cite{lippmann}. However, a well known drawback of the IE is the large cost of numerical integration. In particular, provided the horizontal layering and uniform discretization grid in $x$ and $y$, FFT is applied along $x$ and $y$ and correlation is computed along $z$, so the cost is $O(N_xN_ylog(N_xN_y)N_z^2)$.  Here $N_x$, $N_y$ and $N_z$ are numbers of discretization nodes along $x$, $y$ and $z$, respectively. For non-uniform discretization in $x$ and $y$, discrete Fourier transform is employed and the cost is $O(N_x^2N_y^2N_z^2)$.

Quadratic dependence of the cost on $N_z$ was removed in FDIE approach of \cite{zaslseg, zaslgeo}. FDIE can be viewed as a finite-difference discretization of the volume IE approach. But thanks to the sparsity of the finite-difference discretization, the cost of preconditioner in FDIE is $O(N_xN_ylog(N_xN_y)N_z)$ and $O(N_x^2N_y^2N_z)$ for uniform and non-uniform grids, respectively (similar approach was proposed in \cite{avdeev} for volume IE method). 

For geophysical applications one typically uses fine (not necessarily uniform) grid within the domain of the survey and coarse (exponentially convergent) grid outside (possibly, with imaginary grid steps in PML). Therefore the cost of preconditioners of \cite{zaslseg, zaslgeo, avdeev} is still quadratic with respect to $N_x$ and $N_y$. In this paper we propose a fast algorithm to compute the solution $u$ of Helmholtz equation in layered medium  for arbitrary right hand side part. The algorithm is based on the ideas of cyclic reduction (see \cite{golub, rossi}) which is is closely related to specific method of domain decomposition type. We split the computational grid into connected Cartesian blocks $S^{ij}=S^i_x\times S^j_y\times S_z$ and compute approximate solution $v$ by solving the problem on each block independently with homogeneous Dirichlet boundary conditions. The corresponding residual has support on the interfaces of the blocks $\cup\partial S^{ij}$. Therefore, by performing just local (on $\cup\partial S^{ij}$) computations we obtain Dirichlet conditions $\left(w\right)|_{\partial S^{ij}}$ for $w=u-v$ on each interface $\partial S^{ij}$. The last step of our algorithm is computing $w$ in each $S^{ij}$ under given boundary conditions $\left(w\right)|_{\partial S^{ij}}$. The cost of the algorithm is $O(N_xN_ylog(N_xN_y)N_z)$ and $O(N_x^{3/2}N_y^{3/2}N_z)$ for uniform and non-uniform grids in the domain of the survey, respectively. We also extended our approach for solving Maxwell's equation in the frequency domain. In a very similar way the algorithm can be applied for solving elasticity problems as well.

\section{Scalar problem} We consider the Helmholtz equation in horizontally layered medium
\begin{equation}\label{wave2}
\sigma(z)u_{xx}+\sigma(z)u_{yy}+(\sigma(z)u_z)_z+\lambda u=f\end{equation} in $\R^3$, where $\lambda \in\R$, $\sigma(z)$ is regular enough positive function of $z$.

Let $\left\{h^x_i\right\}^{N_x}_{i=1}$, $\left\{h^y_j\right\}^{N_y}_{j=1}$ and $\left\{h^z_k\right\}^{N_z}_{k=1}$ be (perhaps, complex if PML is considered for $\lambda>0$) grid steps of $S_x$, $S_y$ and $S_z$, respectively. Consider 7-point finite-volume discretization scheme 

\bea \label{fd}
\nonumber\sigma_{k}\left(\frac{u_{i+1,j,k}-u_{i,j,k}}{h^x_i}-\frac{u_{i,j,k}-u_{i-1,j,k}}{h^x_{i-1}}\right)\hh^y_j\hh^z_k+\\
\nonumber\sigma_{k}\left(\frac{u_{i,j+1,k}-u_{i,j,k}}{h^y_j}-\frac{u_{i,j,k}-u_{i,j-1,k}}{h^y_{j-1}}\right)\hh^x_i\hh^z_k+\\
\nonumber\left(\sigma_{k+1/2}\frac{u_{i,j,k+1}-u_{i,j,k}}{h^z_k}-\sigma_{k-1/2}\frac{u_{i,j,k}-u_{i,j,k-1}}{h^z_{k-1}}\right)\hh^x_i\hh^y_j+\\\lambda u_{i,j,k}\hh^x_i\hh^y_j\hh^z_k=f_{i,j,k}\hh^x_i\hh^y_j\hh^z_k,\eea
where $\hh^x_i, \hh^y_j, \hh^z_k$ are lengths of the edges of the control volume $V_{i,j,k}$.

The system (\ref{fd}) can be rewritten as 
\bea\label{fdmtr} \nonumber Au=\\\nonumber \left(A_x\otimes M_y\otimes M_z+M_x\otimes A_y\otimes M_z+M_x\otimes M_y\otimes A_z+\lambda M_x\otimes M_y\otimes M_z\right)u=\\ M_x\otimes M_y\otimes M_zf,\eea
where $A_x$, $A_y$ and $A_z$ are discretizations of $\int_{V_{i,j,k}}{\sigma(z)\frac{\partial^2}{\partial x^2}dV}$, $\int_{V_{i,j,k}}{\sigma(z)\frac{\partial^2}{\partial y^2}dV}$ and $\int_{V_{i,j,k}}{\frac{\partial}{\partial z}\left(\sigma(z){\frac{\partial}{\partial z}}\right)dV}$, respectively, and $M_x$, $M_y$ and $M_z$ are diagonal mass matrices: $M_x=diag(\hh^x_1,\ldots,\hh^x_{N_x})$, for example.

Split the computational grid into connected Cartesian blocks $S^{ij}=S^i_x\times S^j_y\times S_z$ and define $P^{i}_x$ and $P^j_y$ the projection operators to interior nodes of $S^i_x$ and $S^j_y$, respectively. Let $A^i_d$ and $M^i_d$ be projected  operators $A^i_d=P^i_dA_dP^i_d$ and $M^i_d=P^i_dM_dP^i_d$, where $d\in\{x,y,z\}$. In each block $S^{ij}$ we construct the approximate solution $v^{ij}$ satisfying the equation
\bea\label{fdmtrprj} \nonumber A^{ij}v^{ij}=\\\nonumber \left(A^i_x\otimes M^j_y\otimes M_z+M^i_x\otimes A^j_y\otimes M_z+M^i_x\otimes M^j_y\otimes A_z+\lambda M^i_x\otimes M^j_y\otimes M_z\right)v^{ij}=\\ M^i_x\otimes M^j_y\otimes M_zf^{ij}\eea
and homogeneous Dirichlet boundary conditions $v^{ij}|_{\partial S^{ij}}=0$. Here $f^{ij}=P^i_xP^j_yf$. We also define the global approximate solution $v$ as $v|_{S^{ij}}=v^{ij}$ and $v|_{\cup\partial S^{ij}}=0$. 
Denote $\Lambda^i_d$ and $W^i_d$ the diagonal matrix of generalized eigenvalues and matrix of generalized eigenvectors of operator $A^i_d$:
$$A^i_dW^i_d=M^i_dW^i_d\Lambda^i_d, W^i_d|_{\partial S^{i}_d}=0, d\in\{x,y\}.$$
First we Fourier-transform (\ref{fdmtrprj}) with respect to $x$ and $y$:
\bea\label{fdmtrprjfr} \nonumber \left(\Lambda^i_x\otimes I\otimes M_z+I\otimes \Lambda^j_y\otimes M_z+I\otimes I\otimes A_z+\lambda I\otimes I\otimes M_z\right)\tilde{v}^{ij}=\\ I\otimes I\otimes M_z\tilde{f}^{ij},\eea where 
$\tilde{v}^{ij}=(W^i_x)^*\otimes (W^j_y)^*\otimes Iv^{ij}$, $\tilde{f}^{ij}=(W^i_x)^*\otimes (W^j_y)^*\otimes If^{ij}$.
Obviously, the system (\ref{fdmtrprjfr}) can be split (with respect to eigenmode) into multiple tridiagonal systems. For each pair $(i,j)$ and each eigenmode the solution $\tilde{v}^{ij}$ can be computed in $O(N_z)$ operations.

It is easy to see that the residual $r=f-Av$ has a support only on $\cup\partial S^{ij}$. Therefore, thanks to the sparsity of $A$, to obtain $r$ we need $v$ only on $\cup\partial S^{ij}$ as well as at neighboring nodes.
Under known $\tilde{v}^{ij}$, this step can be performed via inverse Fourier transform.

Obviously, the difference $w=u-v$ satisfies the homogeneous equation $A^{ij}w=0$ in each $S^{ij}$ with some Dirichlet boundary conditions $w|_{\partial S^{ij}}$. Thanks to the sparsity of $r$, the latter can be obtained in $O((N^2_x+N^2_y)N_z)$ operations via Fourier transform of (\ref{fdmtr}). In fact, let $\Lambda_d$ and $W_d$ be diagonal matrix of generalized eigenvalues and matrix of generalized eigenvectors of operator $A_d$, respectively: $A_dW_d=M_dW_d\Lambda_d, d\in\{x,y\},~W_d|_{\partial S}=0$. Then $\tilde{w}=(W_x)^*\otimes (W_y)^*\otimes I v$ satisfies 
\bea\label{fdmtrprjfrerr} \nonumber \left(\Lambda_x\otimes I\otimes M_z+I\otimes \Lambda_y\otimes M_z+I\otimes I\otimes A_z+\lambda I\otimes I\otimes M_z\right)\tilde{w}=\\ I\otimes I\otimes M_z\tilde{r},\eea where Fourier transform $\tilde{r}=(W_x)^*\otimes (W_y)^*\otimes I r$ can be performed in $O((N^2_x+N^2_y)N_z)$ operations thanks to the sparsity of $r$. The system (\ref{fdmtrprjfrerr}) can be split into $N_xN_y$ tridiagonal systems and each of them can be solved in $O(N_z)$ operations. Finally, to obtain $u|_{\cup\partial S^{ij}}=w|_{\cup\partial S^{ij}}=\left(W_x\otimes W_y\otimes I \tilde{w}\right)|_{\cup\partial S^{ij}}$ we need to perform another $O((N^2_x+N^2_y)N_z)$ operations of inverse Fourier transform on the $\cup\partial S^{ij}$. 

To compute $u=v+w$ in the interior nodes of each $S^{ij}$, we first solve the Fourier transformed equation 
\bea\label{fdmtrprjfrerrloc} \nonumber \left(\Lambda^i_x\otimes I\otimes M_z+I\otimes \Lambda^j_y\otimes M_z+I\otimes I\otimes A_z+\lambda I\otimes I\otimes M_z\right)\tilde{w}^{ij}=0\eea
in $S^{ij}$, where $\tilde{w}^{ij}=(W^i_x)^*\otimes (W^j_y)^*\otimes I P^i_xP^j_yw$ satisfies the computed boundary conditions $w|_{\partial S^{ij}}$. Finally, in each $S^{ij}$ we perform inverse Fourier transform and obtain $P^i_xP^j_yu=W^i_x\otimes W^j_y\otimes I(\tilde{u}^{ij}+\tilde{w}^{ij})$.

Our algorithm can be summarized as follows

\begin{enumerate}
\item for each $i,j$ perform Fourier transform and compute $\tilde{f}^{ij}=(W^i_x)^*\otimes (W^j_y)^*\otimes IP^i_xP^j_yf$ 
\item for each $i,j$ solve the problem (\ref{fdmtrprjfr}) and compute $\tilde{v}^{ij}$
\item via inverse Fourier transform, compute $v$ at the nodes adjacent to $\cup\partial S^{ij}$ .
\item compute residual $r=f-Av$ at $\cup\partial S^{ij}$
\item compute $w|_{\cup\partial S^{ij}}$ via Fourier transform with respect to $W_x$ and $W_y$
\item under known $w|_{\cup\partial S^{ij}}$, compute $\tilde{w}^{ij}$ for each $i,j$
\item compute $u=v+w$ via inverse Fourier transform of $\tilde{v}^{ij}+\tilde{w}^{ij}$ in each $S^{ij}$
\end{enumerate}

The computational costs of the second, the fourth and the fifth steps are $O(N_xN_yN_z)$, $O((N_x+N_y)N_z)$ and $O((N^2_x+N^2_y)N_z)$, respectively. The costs of the remaining steps of the algorithm depend on how we split our grid $S$ into blocks $S^{ij}$. For the case of uniform grids in the domain of the survey we split each of $S_d, d\in\{x,y\}$ into three subgrids: $S^1_d$ and $S^3_d$ correspond to the non-uniform exponentially convergent grids and $S^2_d$ covers the uniform part. To maintain the same approximation error within the domain of the survey and outside, further we assume that the number of nodes in exponentially convergent grids along $x$ and $y$ is proportional to logarithm of uniform parts of $S_x$ and $S_y$, respectively. For the case of general non-uniform grids we split $S_d$ ($d\in\{x,y\}$) into $\sqrt{N_d}$ subgrids with $\sqrt{N_d}$ nodes in each of them. For the first type of gridding we compute Fourier transforms in uniform parts via FFT. Therefore, the costs of the first, the sixth and the seventh steps are $O(N_xN_ylog(N_xN_y)N_z)$. For general non-uniform gridding, the cost of Fourier transforms on each $S^{ij}$ is $O((\sqrt{N_x})^2(\sqrt{N_y})^2N_z)$ and the total number of blocks is $\sqrt{N_x}\sqrt{N_y}$. Hence, the overall costs of the first, the sixth and the seventh steps are $O(N^{3/2}_xN^{3/2}_yN_z)$. The cost of the third step is  $O(N_xN_yN_z)$ and $O((N^{3/2}_x+N^{3/2}_y)N_z)$ for two type of griddings above.

\begin{remark}
By extending our approach to multi-level cyclic reduction (i.e. by employing multiple mutually embedded partitions of $S$), we can reduce the asymptotics of the cost from $N^{3/2}$ to $Nlog(N)$ along each direction. But the constant factor is significantly greater due to multiple (on each level of embedding) solution of block-tridiagonal systems.
\end{remark}

\section{Maxwell equations}

Consider the Maxwell equation for magnetic field in the layered medium:
\begin{eqnarray}
\label{Maxwell1} \nabla \times\E &=& {\rm i}\omega\mu\H\,,
\nonumber \\
\nabla\times  {\H} &=& \rho^{-1}\E+\J'
\end{eqnarray}

Here we made an assumption of negligible displacement current which is typical for large-scale problems of geophysical prospecting. 

Following \cite{davyd}, we
discretize Maxwell equations using Lebedev grid \cite{lebedev},
which is the counterpart of the Yee grid \cite{yee} for anisotropic problems.

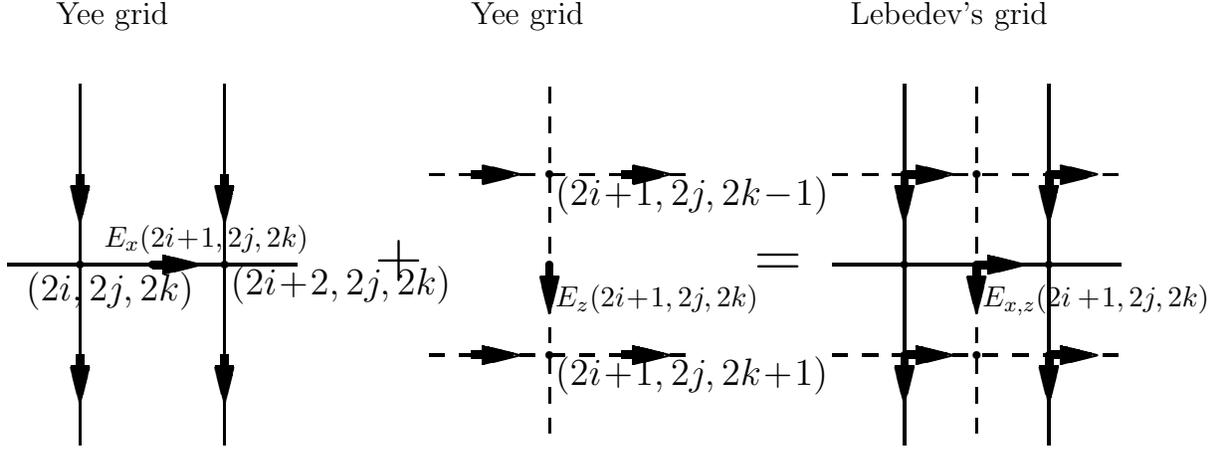
\begin{figure}
\unitlength=0.8mm
\begin{picture}(190,80)(0,-40)
% ris. 1
\thicklines
\put(0,0){\line(1,0){48}}
%===============begin arrow
\put(90,0){\vdarrow}
\put(12,15){\vdarrow}
\put(36,15){\vdarrow}
\put(12,-15){\vdarrow}
\put(36,-15){\vdarrow}
\put(149,15){\vdarrow}
\put(173,15){\vdarrow}
\put(149,-15){\vdarrow}
\put(173,-15){\vdarrow}
%=============end arrow
%===============begin arrow
\put(24,0){\hrarrow}
%=============end arrow
\put(0,0){\line(1,0){48}}
%===============begin horiz. arrows
\put(77,15){\hrarrow}
\put(102,15){\hrarrow}
\put(77,-15){\hrarrow}
\put(102,-15){\hrarrow}
\put(149,15){\hrarrow}
\put(173,15){\hrarrow}
\put(149,-15){\hrarrow}
\put(173,-15){\hrarrow}
%=============end horiz. arrow
%===============begin arrow
\put(161,0){\hrarrow}
\put(161,0){\vdarrow}
%=============end arrow
%
\put(12,-30){\line(0,1){60}}
\put(36,-30){\line(0,1){60}}
\put(12,0){\circle*{1.2}}
\put(24,0){\circle*{1.2}}
\put(36,0){\circle*{1.2}}
\put(16,3){\mbox{\small{$E_x(2i\!+\!1,2j,2k)$}}}
\put(37,-5){\mbox{\large{${(2i\!+\!2,2j,2k)}$}}}
\put(3,-6){\mbox{\large{$(2i,2j,2k)$}}}
\put(8,40){\mbox{Yee grid}}
\put(61,-2){\mbox{\Huge $+$}}
% ris. 2
\multiput(90,-28)(0,5){12}{\line(0,1){2.5}}
\multiput(70,-15)(5,0){9}{\line(1,0){2.5}}
\multiput(70,15)(5,0){9}{\line(1,0){2.5}}
\put(90,15){\circle*{1.2}}
\put(90,0){\circle*{1.2}}
\put(90,-15){\circle*{1.2}}
\put(91,-7){\mbox{\small{$E_z(2i\!+\!1,2j,2k)$}}}
\put(91,10){\mbox{\large{$(2i\!+\!1,2j,2k\!-\!1)$}}}
\put(91,-20){\mbox{\large{$(2i\!+\!1,2j,2k\!+\!1)$}}}
\put(77,40){\mbox{Yee grid}}
\put(124,-2){\mbox{\Huge $=$}}
% ris.3
\put(137,0){\line(1,0){48}}
\put(149,-30){\line(0,1){60}}
\put(173,-30){\line(0,1){60}}
\put(149,0){\circle*{1.2}}
\put(149,15){\circle*{1.2}}
\put(149,-15){\circle*{1.2}}
\put(173,0){\circle*{1.2}}
\put(173,15){\circle*{1.2}}
\put(173,-15){\circle*{1.2}}
\multiput(161,-28)(0,5){12}{\line(0,1){2.5}}
\multiput(137,-15)(5,0){10}{\line(1,0){2.5}}
\multiput(137,15)(5,0){10}{\line(1,0){2.5}}
\put(161,15){\circle*{1.2}}
\put(161,0){\circle*{1.2}}
\put(161,-15){\circle*{1.2}}
%\put(162,3){\mbox{\small{$E_x(2i+\!1,2j,2k)$}}}
\put(162,-7){\mbox{\small{$E_{x,z}(2i+\!1,2j,2k)$}}}
\put(140,40){\mbox{Lebedev's grid}}
\end{picture}
\caption{2-D cross-section (in the plane $Oxz$, for even nodes $y$)
of Lebedev's staggered grid consisting of two clusters.
The crossings of the lines of the same type form the subgrid
$P$, whereas the crossings of the lines of different types
form the subgrid $R$.}
\label{lebyee}
\end{figure}

We consider a bounded computational domain
$\Omega=[x_{1},x_{N_x}]\times[y_{1},y_{N_y}]\times[z_{1},z_{N_z}]$
and introduce a Cartesian three-dimensional (3D) grid as follows:
\begin{eqnarray}
M&=&\left(\r_{i,j,k}\right)\,,\quad
i=1,\cdots,N_x\nonumber\\
&&j=1,\cdots,N_y,\,,\quad k=1,\cdots,N_z\,,
\end{eqnarray}
%$S=\{(ih_1;jh_2;kh_3)\}^{N_x,N_y,N_z}_{i,j,k=0}$.
where $\r_{i,j,k}=(x_i,y_j,z_k)$, $N_x$, $N_y$ and $N_z$ are the
number of grid nodes in the $x-$, $y-$ and $z-$directions. The
Lebedev $P$-grid is defined as a sub-grid of $M$ with even sum of
indices $(i+j+k)$ and the Lebedev $R$-grid contains the remaining
nodes of $M$. For example, for even $k$ the grid in the $xy-$plane
is shown in Fig.~\ref{lebyee}.

In other words, each of the Lebedev $P$- and $R$-grids is a
combination of two Yee grids in the two-dimensional (2D)
configuration and four Yee grids in 3D configuration (see Fig. \ref{lebyee}). For any
one-dimensional (1D) grid $M^1$ we define primary and dual grids as
sub-grids of $M^1$ with even and odd indices, respectively.
Therefore, in the $xy-$plane and for even k, the Lebedev $P$-grid is
a superposition of the tensor product of two primary grids as well
as the tensor product of two dual grids. On the other hand, for odd
$k$, the Lebedev $P$-grid is a superposition of the tensor product
of a primary (along $x$) and a dual (along $y$) grid as well as the
tensor product of a dual (along $x$) and a primary (along $y$) grid.
In fact, according to the definition of the $P$-grid and for even
$k$, it consists of vertices with either both even $i$ and $j$
(which is a tensor product of two primary grids) or both odd $i$ and
$j$ (which is a tensor product of two dual grids). Similarly, for
odd $k$ we have a tensor product of primary and dual grids as well
as a tensor product of dual and primary grids.

Let all three components of the finite-difference vector magnetic
field $\H_h=(H_h^x,H_h^y,H_h^z)$ be defined at the same nodes of the
$P$-grid, and similarly, all three components of the
finite-difference vector electric field $\E_h=(E_h^x,E_h^y,E_h^z)$
be defined at the same nodes of the $R$-grid. Hence there is no need
to perform any interpolation with our scheme to handle anisotropic
Ohm's law. We define finite-difference derivatives on both grids
along the $x-$direction as follows:
\begin{equation}\label{fin}
(f_x)_{i,j,k}=\frac{f_{i+1,j,k}-f_{i-1,j,k}}{x_{i+1}-x_{i-1}}\,,\\
\end{equation}
and similarly along $y-$ and $z-$direction. Note that equation
(\ref{fin}) performs mapping from $P$-grid to $R$-grid and vice
versa. This allows us to obtain a self-consistent discrete system of
Maxwell equations
\begin{eqnarray} \label{Maxwell1fd} \nabla_h
\times\E &=& {\rm i}\omega\mu\H\,,
\nonumber \\
\nabla_h\times  {\H} &=& \rho\E+\J',
\end{eqnarray}
here we use the same notations for discrete fields as for continuous ones.

 After
eliminating $\E$ from  (\ref{Maxwell1fd}) we obtain the linear
system of equations for the magnetic field vector $\H$
\begin{equation}\label{abstr}\TA_h\H=\f\,,
\end{equation}
where
\[\TA_h{\H}=\nabla_h\times\rho\nabla_h\times {\H}-{\rm
i}\omega \mu\H\]
and
\begin{equation}\label{rhs}
{\f}=\nabla_h\times\rho  {\J'}\,. \end{equation}

The boundary of the computational domain $\partial \Omega$ contains
nodes of both $P$- and $R$-grids. We replace the boundary conditions
at infinity by:
\begin{equation}
\H\times\n|_{\partial \Omega\cap P}=0, \quad \E\times\n|_{\partial
\Omega \cap R}=0\,.
\end{equation}
As pointed out in \cite{davyd}, this allows us to
decrease the computational domain by a factor of 2 (compared to the
conventional boundary conditions $\H|_{\partial\Omega}=0$ or
$\E|_{\partial\Omega}=0$) without losing accuracy.

Consider the Fourier transform (with respect to $x$ and $y$) of equations (\ref{Maxwell1fd}). Due to presence of mixed derivatives in (\ref{Maxwell1fd}), the equations in Fourier domain can not be represented in such a simple Krononecker product form as (\ref{fdmtrprjfrerr}). According to the definition of the
Lebedev grid, for even $z$ grid node $k$ in the $xy$-plane, the
magnetic field is defined on a tensor product of two primary grids
and a tensor product of two dual grids. Denote them by ${\H}^{pp}$
and ${\H}^{dd}$, where the first and second indices mean the type of
the grid along $x$ and $y$, respectively. Similarly, ${\H}^{pd}$ and
${\H}^{dp}$ are defined for odd $k$. Define eigenvalues and
eigenfunctions of 1D operators on each grid with Dirichlet
conditions on primary grid and Neumann conditions on dual grid:
 $$- \phi^{\alpha,l}_{xx}=(\lambda^l)^2\phi^{\alpha,l}, \quad - \psi^{\beta,m}_{yy}=(\nu^m)^2\psi^{\beta,m},$$
where $\alpha, \beta\in\{p,d\}$. One can derive the relations between
primary and dual eigenfunctions: $\phi^{d,l}=\lambda^l\phi^{p,l}_x$
and $\psi^{d,m}=\nu^m\psi^{p,m}_y$. For each $k$, we expand $\H$ and
$\f$ as a sum
$${\H}^{\alpha\beta}_k=\sum_{lm}{\tilde{{\H}}^{\alpha\beta,lm}_k\phi^{\alpha,l}\psi^{\beta,m}},$$
$${\f}^{\alpha\beta}=\sum_{lm}{\tilde{{\f}}^{\alpha\beta,lm}\phi^{\alpha,l}\psi^{\beta,m}}.$$
Here $\alpha\beta\in\{pp,dd\}$ for even $k$ and
$\alpha\beta\in\{pd,dp\}$ for odd $k$. In this case, the equations
for the expansion coefficients have the following form (we omit the
indices $l$ and $m$ of $\lambda$, $\nu$, $\tilde{{\f}}$ and
$\tilde{\H}$ to be more concise):

For even $k$ and for the primary-primary grid:

$$
\begin{array}{c}
c^j_3\left(\lambda\nu\tilde{H}^{dd}_{k,y}+\nu^2\tilde{H}^{pp}_{k,x}\right)-\\
-\frac{1}{h_k}\left(c^{k+1}_2\frac{\tilde{H}^{pp}_{k+2,x}-\tilde{H}^{pp}_{k,x}}{h_{k+1}}-
c^{k-1}_2\frac{\tilde{H}^{pp}_{k,x}-\tilde{H}^{pp}_{k-2,x}}{h_{k-1}}\right)-\\
-\lambda
\frac{c^{k+1}_2\tilde{H}^{dp}_{k+1,z}-c^{k-1}_2\tilde{H}^{dp}_{k-1,z}}{h_k}-i\omega\mu\tilde{H}^{pp}_{k,x}=\tilde{f}^{pp}_x\\
c^j_3\left(\lambda\nu\tilde{H}^{dd}_{k,x}+\lambda^2\tilde{H}^{pp}_{k,y}\right)-\\
-\frac{1}{h_k}\left(c^{k+1}_1\frac{\tilde{H}^{pp}_{k+2,y}-\tilde{H}^{pp}_{k,y}}{h_{k+1}}-
c^{k-1}_1\frac{\tilde{H}^{pp}_{k,y}-\tilde{H}^{pp}_{k-2,y}}{h_{k-1}}\right)-\\
-\nu
\frac{c^{k+1}_1\tilde{H}^{pd}_{k+1,z}-c^{k-1}_1\tilde{H}^{pd}_{k-1,z}}{h_k}-i\omega\mu\tilde{H}^{pp}_{k,y}=\tilde{f}^{pp}_y\\
\left(c^j_2\lambda^2+c^j_1\nu^2\right)\tilde{H}^{pp}_{k,z}-\lambda
c^j_2 \frac{\tilde{H}^{dp}_{k+1,x}-\tilde{H}^{dp}_{k-1,x}}{h_k}-\\
-\nu c^j_1\frac{\tilde{H}^{pd}_{k+1,y}-\tilde{H}^{pd}_{k-1,y}}{h_k}
-i\omega\mu\tilde{H}^{pp}_{k,z}=\tilde{f}^{pp}_z.
\end{array}
$$

Here and below $c^{k}_i=c_i(z_k)$ and finite-difference derivatives are defined by (\ref{fin}). For even $k$ and for the dual-dual grid:
$$
\begin{array}{c}
c^j_3\left(\lambda\nu\tilde{H}^{pp}_{k,y}+\nu^2\tilde{H}^{dd}_{k,x}\right)-\\
-\frac{1}{h_k}\left(c^{k+1}_2\frac{\tilde{H}^{dd}_{k+2,x}-\tilde{H}^{dd}_{k,x}}{h_{k+1}}-
c^{k-1}_2\frac{\tilde{H}^{dd}_{k,x}-\tilde{H}^{dd}_{k-2,x}}{h_{k-1}}\right)+\\
+\lambda
\frac{c^{k+1}_2\tilde{H}^{dp}_{k+1,z}-c^{k-1}_2\tilde{H}^{dp}_{k-1,z}}{h_k}-i\omega\mu\tilde{H}^{dd}_{k,x}=\tilde{f}^{dd}_x\\
c^j_3\left(\lambda\nu\tilde{H}^{pp}_{k,x}+\lambda^2\tilde{H}^{dd}_{k,y}\right)-\\
-\frac{1}{h_k}\left(c^{k+1}_1\frac{\tilde{H}^{dd}_{k+2,y}-\tilde{H}^{dd}_{k,y}}{h_{k+1}}-
c^{k-1}_1\frac{\tilde{H}^{dd}_{k,y}-\tilde{H}^{dd}_{k-2,y}}{h_{k-1}}\right)+\\
+\nu
\frac{c^{k+1}_1\tilde{H}^{pd}_{k+1,z}-c^{k-1}_1\tilde{H}^{pd}_{k-1,z}}{h_k}-i\omega\mu\tilde{H}^{dd}_{k,y}=\tilde{f}^{dd}_y\\
\left(c^j_2\lambda^2+c^j_1\nu^2\right)\tilde{H}^{dd}_{k,z}+\lambda
c^j_2 \frac{\tilde{H}^{dp}_{k+1,x}-\tilde{H}^{dp}_{k-1,x}}{h_k}+\\
+\nu c^j_1\frac{\tilde{H}^{pd}_{k+1,y}-\tilde{H}^{pd}_{k-1,y}}{h_k}
-i\omega\mu\tilde{H}^{dd}_{k,z}=\tilde{f}^{dd}_z.
\end{array}
$$

For odd $k$ and for the primary-dual grid:
$$
\begin{array}{c}
c^j_3\left(-\lambda\nu\tilde{H}^{dp}_{k,y}+\nu^2\tilde{H}^{pd}_{k,x}\right)-\\
-\frac{1}{h_k}\left(c^{k+1}_2\frac{\tilde{H}^{pd}_{k+2,x}-\tilde{H}^{pd}_{k,x}}{h_{k+1}}-
c^{k-1}_2\frac{\tilde{H}^{pd}_{k,x}-\tilde{H}^{pd}_{k-2,x}}{h_{k-1}}\right)-\\
-\lambda
\frac{c^{k+1}_2\tilde{H}^{dd}_{k+1,z}-c^{k-1}_2\tilde{H}^{dd}_{k-1,z}}{h_k}-i\omega\mu\tilde{H}^{pd}_{k,x}=\tilde{f}^{pd}_x\\
c^j_3\left(-\lambda\nu\tilde{H}^{dp}_{k,x}+\lambda^2\tilde{H}^{pd}_{k,y}\right)-\\
-\frac{1}{h_k}\left(c^{k+1}_1\frac{\tilde{H}^{pd}_{k+2,y}-\tilde{H}^{pd}_{k,y}}{h_{k+1}}-
c^{k-1}_1\frac{\tilde{H}^{pd}_{k,y}-\tilde{H}^{pd}_{k-2,y}}{h_{k-1}}\right)+\\
+\nu
\frac{c^{k+1}_1\tilde{H}^{pp}_{k+1,z}-c^{k-1}_1\tilde{H}^{pp}_{k-1,z}}{h_k}-i\omega\mu\tilde{H}^{pd}_{k,y}=\tilde{f}^{pd}_y\\
\left(c^j_2\lambda^2+c^j_1\nu^2\right)\tilde{H}^{pd}_{k,z}-\lambda
c^j_2 \frac{\tilde{H}^{dd}_{k+1,x}-\tilde{H}^{dd}_{k-1,x}}{h_k}+\\
+\nu c^j_1\frac{\tilde{H}^{pp}_{k+1,y}-\tilde{H}^{pp}_{k-1,y}}{h_k}
-i\omega\mu\tilde{H}^{pd}_{k,z}=\tilde{f}^{pd}_z.
\end{array}
$$

Finally for odd $k$ and for the dual-primary grid:
$$
\begin{array}{c}
c^j_3\left(-\lambda\nu\tilde{H}^{pd}_{k,y}+\nu^2\tilde{H}^{dp}_{k,x}\right)-\\
-\frac{1}{h_k}\left(c^{k+1}_2\frac{\tilde{H}^{dp}_{k+2,x}-\tilde{H}^{dp}_{k,x}}{h_{k+1}}-
c^{k-1}_2\frac{\tilde{H}^{dp}_{k,x}-\tilde{H}^{dp}_{k-2,x}}{h_{k-1}}\right)+\\
+\lambda
\frac{c^{k+1}_2\tilde{H}^{pp}_{k+1,z}-c^{k-1}_2\tilde{H}^{pp}_{k-1,z}}{h_k}-i\omega\mu\tilde{H}^{dp}_{k,x}=\tilde{f}^{dp}_x\\
c^j_3\left(-\lambda\nu\tilde{H}^{pd}_{k,x}+\lambda^2\tilde{H}^{dp}_{k,y}\right)-\\
-\frac{1}{h_k}\left(c^{k+1}_1\frac{\tilde{H}^{dp}_{k+2,y}-\tilde{H}^{dp}_{k,y}}{h_{k+1}}-
c^{k-1}_1\frac{\tilde{H}^{dp}_{k,y}-\tilde{H}^{dp}_{k-2,y}}{h_{k-1}}\right)-\\
-\nu
\frac{c^{k+1}_1\tilde{H}^{dd}_{k+1,z}-c^{k-1}_1\tilde{H}^{dd}_{k-1,z}}{h_k}-i\omega\mu\tilde{H}^{dp}_{k,y}=\tilde{f}^{dp}_y\\
\left(c^j_2\lambda^2+c^j_1\nu^2\right)\tilde{H}^{dp}_{k,z}+\lambda
c^j_2 \frac{\tilde{H}^{pp}_{k+1,x}-\tilde{H}^{pp}_{k-1,x}}{h_k}-\\
-\nu c^j_1\frac{\tilde{H}^{dd}_{k+1,y}-\tilde{H}^{dd}_{k-1,y}}{h_k}
-i\omega\mu\tilde{H}^{dp}_{k,z}=\tilde{f}^{dp}_z.
\end{array}
$$

Thus for a particular harmonic $lm$ we obtain a block 5-diagonal
linear system where each block is 6x6 matrix. In the case of a
diagonal $\rho$, the system decouples into two block 5-diagonal
linear systems with 3x3 blocks. Further simplifications can be made:
each $\tilde{H}_{k,z}$ may be expressed in terms of
$\tilde{H}_{k,x}$ and $\tilde{H}_{k,y}$ in accordance with the third
equation of the above set of equations. Substituting for
$\tilde{H}^{\alpha\beta}_{k,z}$ by that expression in the remaining
equations we obtain a block 3-diagonal linear system with 2x2
blocks.

Similar to the scalar Helmholtz problem, we split the solution $\H$ (as well as $\E$) into two parts: $\H=\H^1+\H^2$ ($\E=\E^1+\E^2$), where $\H^1$ satisfies (\ref{Maxwell1fd}) in each $S^{ij}$ and the following homogeneous boundary conditions
\begin{equation}
\H^1\times\n|_{\partial S^{ij}\cap P}=0, \quad \E^1\times\n|_{\partial
S^{ij} \cap R}=0\,.
\end{equation}
Due to the coupling of $P-$ and $R$-grids, in contrast with the Helmholtz problem, the boundary conditions are imposed on both these grids. Consequently, to obtain $\H^2$, the boundary conditions 
\begin{equation}
\H^2\times\n|_{\partial S^{ij}\cap P}, \quad \E^2\times\n|_{\partial
S^{ij} \cap R}\,.
\end{equation}
are supposed to be computed on the fifth step of our two-level cyclic reduction algorithm (see scalar case). Other steps remain the same.

\end{document}